\documentclass[a4paper,10pt]{amsart}						

	\usepackage{amsthm, amsfonts, verbatim,appendix}
	\usepackage{amsmath}   								
	\usepackage{mathrsfs}	
	\usepackage{graphicx}		

	\newtheorem{thm}{Theorem}

	\theoremstyle{definition}

  	\newtheorem{rem}[thm]{Remark}

	\newcommand{\M}{\mathcal{M}}
	\newcommand{\Mbar}{\overline{\mathcal{M}}}

	\newcommand{\g}{\mathfrak{g}}
	
	\newcommand{\D}{\mathfrak{D}}

	\newcommand{\cc}{\mathfrak{c}}

	\newcommand{\DD}{\mathcal{D}}
	\newcommand{\GG}{\mathcal{G}}	
	\newcommand{\YY}{\mathcal{Y}}	
	\newcommand{\Zz}{\mathcal{Z}}
	
	\newcommand{\EE}{\mathscr{E}}
	\newcommand{\II}{\mathscr{I}}
	\newcommand{\Li}{\mathscr{L}}
	\newcommand{\MM}{\mathscr{M}}
	\newcommand{\OO}{\mathscr{O}}
	\newcommand{\VV}{\mathscr{V}}

	\newcommand{\PP}{\mathbb{P}}
	\newcommand{\ZZ}{\mathbb{Z}}
	\newcommand{\QQ}{\mathbb{Q}}

	\title{Double points of plane models in $\Mbar_{g,1}$}
	\author{Nicola Tarasca}
	
	\address{Humboldt Universit\"at zu Berlin and Berlin Mathematical School}
	\email{tarasca@math.hu-berlin.de}

\begin{document}

\begin{abstract}
The aim of this paper is to compute the class of the closure of the effective divisor $\D^2_6$ in $\M_{6,1}$ given by pointed curves $[C,p]$ with a sextic plane model mapping $p$ to a double point.
Such a divisor generates an extremal ray in the pseudoeffective cone of $\Mbar_{6,1}$ as shown by Jensen.
A general result on some families of linear series with adjusted Brill-Noether number $0$ or $-1$ is introduced to complete the computation.
\end{abstract}

\maketitle

The birational geometry of an algebraic variety is encoded in its cone of effective divisors. Nowadays a major problem is to determine the effective cone of moduli spaces of curves.

Let $\mathcal{GP}^1_4$ be the Gieseker-Petri divisor in $\M_6$ given by curves with a $\g^1_4$ violating the Petri condition. The class
\[
\left[ \overline{\mathcal{GP}}^1_4 \right] =94 \lambda -12 \delta_0 -50 \delta_1-78 \delta_2- 88\delta_3 \in {\rm Pic}_\QQ(\Mbar_6)
\]
is computed in \cite{MR910206} where classes of Brill-Noether divisors and Gieseker-Petri divisors are determined for arbitrary genera in order to prove the general type of $\Mbar_g$ for $g\geq 24$. 

Now let $\D^2_d$ be the divisor in $\M_{g,1}$ defined as the locus of smooth pointed curves $[C,p]$ with a net $\g^2_d$ of Brill-Noether number $0$ mapping $p$ to a double point. That is
\[
\D^2_{d} := \left\{ [C,p]\in \M_{g,1} \, | \, \exists \, l \in G^2_d(C) \, \mbox{with} \,\, l(-p-x)\in G^1_{d-2}(C) \,\mbox{where} \,  x\in  C, x\not= p  \,  \right\}
\]
for values of $g,d$ such that $g=3(g-d+2)$. Recently Jensen has shown that $\overline{\D}^2_6$ and the pull-back  of $\overline{\mathcal{GP}}^1_4$ to $\Mbar_{6,1}$ generate extremal rays of the pseudoeffective cone of $\Mbar_{6,1}$ (see \cite{Jensen}). Our aim is to prove the following theorem.
\begin{thm}
\label{main}
The class of the divisor $\overline{\D}^2_6 \subset \Mbar_{6,1}$ is
\[
\left[ \overline{\D}^2_6 \right] = 62 \lambda +4 \psi - 8\delta_0 - 30\delta_1 - 52\delta_2 -60 \delta_ 3- 54 \delta_ 4- 34\delta_5 \in {\rm Pic}_\QQ(\Mbar_{6,1}).
\]
\end{thm}

A mix of a Porteous-type argument, the method of test curves and a pull-back to rational pointed curves will lead to the result. Following a method described in \cite{khosla}, we realize $ \overline{\D}^2_d$ in $\M^{\rm irr}_{g,1}$ as the push-forward of a degeneracy locus of a map of vector bundles over $\GG^2_d(\M^{\rm irr}_{g,1})$. This will give us the coefficients of $\lambda$, $\psi$ and $\delta_0$ for the class of  $\overline{\D}^2_d$ in general. Intersecting $\overline{\D}^2_d$ with carefully chosen one-dimensional families of curves will produce relations to determine the coefficients of $\delta_1$ and $\delta_{g-1}$. Finally in the case $g=6$ we will get enough relations to find the other coefficients by pulling-back to the moduli space of stable pointed rational curves in the spirit of \cite[\S 3]{MR910206}.

To complete our computation we obtain a general result on some families of linear series on pointed curves with adjusted Brill-Noether number $\rho=0$ that morally excludes further ramifications on such families. 

\vskip2pc %%%%%

\begin{thm}
\label{ramif0}
Let $(C,y)$ be a general pointed curve of genus $g>1$. Let $l$ be a $\g^r_d$ on $C$ with $r\geq 2$ and adjusted Brill-Noether number $\rho(C,y)=0$. Denote by $(a_0, a_1, \dots, a_r)$ the vanishing sequence of $l$ at $y$. Then $l(-a_{i} y)$ is base-point free for $i=0,\dots,r-1$.
\end{thm}

\begin{comment}
\begin{thm}
\label{ramif00}
Let $(C,y)$ be a general pointed curve of genus $g>1$. Let $l\in G^r_d(C)$ be a linear series on $C$ of dimension $r \geq 2$ and $(a_0, \dots, a_r)$ be the vanishing sequence of $l$ at $y$. If 
$
g-(r+1)(g-d+r)=\sum_{i=0}^r (a_i-i)
$
then $l(-a_{i} y)$ is base-point free for $i=0,\dots,r-1$.
\end{thm}
\end{comment}

\noindent For instance if $C$ is a general curve of genus $4$ and $l\in G^2_5(C)$ has vanishing sequence $(0,1,3)$ at a general point $p$ in $C$, then $l(-p)$ is base-point free. 

Using the irreducibility of the families of linear series with adjusted Brill-Noether number $-1$ (\cite{MR985853}), we get a similar statement for an arbitrary point on the general curve in such families.

\begin{thm}
\label{bpf-1}
Let $C$ be a general curve of genus $g>2$. Let $l$ be a $\g^r_d$ on $C$ with $r\geq 2$ and adjusted Brill-Noether number $\rho(C,y)=-1$ at an arbitrary point $y$. Denote by $(a_0, a_1, \dots, a_r)$ the vanishing sequence of $l$ at $y$. Then $l(-a_1 y)$ is base-point free.
\end{thm}

As a verification of Thm. \ref{main}, let us note that the class of $\overline{\D}^2_6$ is not a linear combination of the class of the Gieseker-Petri divisor $\mathcal{GP}^1_4$ and the class of the divisor $\mathcal{W}$ of Weierstrass points computed in \cite{MR1016424}
\[
[\mathcal{W}] = - \lambda + 21 \psi - 15\delta_1 - 10\delta_2 - 6\delta_3 - 3\delta_4 -\delta_5  \in {\rm Pic}_\QQ(\Mbar_{6,1}).
\]

After briefly recalling in the next section some basic results about limit linear series and enumerative geometry on the general curve, we prove Thm. \ref{ramif0} and Thm. \ref{bpf-1} in section \ref{secram}. Finally in section \ref{class} we prove a general version of Thm. \ref{main}.

\vskip1pc

\noindent {\bf Acknowledgment} This work is part of my PhD thesis. I am grateful to my advisor Gavril Farkas for his guidance. I have been supported by the Graduierten\-kolleg 870 and the Berlin Mathematical School.

\section{Limit linear series and enumerative geometry}

We use throughout Eisenbud and Harris's theory of limit linear series (see \cite{MR846932}). Let us recall some basic definitions and results.

\subsection{Linear series on pointed curves}

Let $C$ be a complex smooth projective curve of genus $g$ and $l=(\Li, V)$ a linear series of type $\g^r_d$ on $C$, that is $\Li \in {\rm Pic}^d(C)$ and $V\subset H^0(\Li)$ is a subspace of vector-space dimension $r+1$. The {\it vanishing sequence} $a^l(p): 0\leq a_0 <\dots <a_r\leq d$ of $l$ at a point $p\in C$ is defined as the sequence of distinct order of vanishing of sections in $V$ at $p$, and the {\it ramification sequence} $\alpha^l(p): 0\leq \alpha_0 \leq \dots \leq \alpha_r\leq d-r$ as $\alpha_i := a_i-i$, for $i=0,\dots,r$. The {\it weight} $w^l(p)$ will be the sum of the $\alpha_i$'s.

Given an $n$-pointed curve $(C,p_1,\dots,p_n)$ of genus $g$ and $l$ a $\g^r_d$ on $C$, the {\it adjusted Brill-Noether number} is 
\[
\rho(C,p_1,\dots p_n) = \rho(g,r,d,\alpha^l(p_1),\dots,\alpha^l(p_n)) := g-(r+1)(g-d+r)-\sum_{i,j} \alpha^l_j(p_i).
\]

\subsection{Counting linear series on the general curve}

Let $C$ be a general curve of genus $g$ and consider $r,d$ such that $\rho(g,r,d)=0$. Then by Brill-Noether theory, the curve $C$ admits only a finite number of $\g^r_d$'s computed by the {\it Castelnuovo number}
\[
N_{g,r,d} := g!\prod_{i=0}^r \frac{i!}{(g-d+r+i)!}.
\]
Furthermore let $(C,p)$ be a general pointed curve of genus $g$ and let $\overline{\alpha}=(\alpha_0, \dots,\alpha_r)$ be a Schubert index of type $r,d$ (that is $0\leq \alpha_0 \leq \dots \leq \alpha_r \leq d-r$) such that $\rho(g,r,d, \overline{\alpha})=0$. Then by \cite[Prop. 1.2]{MR910206}, the curve $C$ admits a $\g^r_d$ with ramification sequence $\overline{\alpha}$ at the point $p$ if and only if $\alpha_0 +g-d+r\geq 0$. When such linear series exist, there is a finite number of them counted by the following formula
\[
N_{g,r,d,\overline{\alpha}} := g! \frac{\prod_{i<j} (\alpha_j - \alpha_i +j-i)}{\prod_{i=0}^r (g-d+r+\alpha_i+i)!}.
\]

\subsection{Limit linear series}

For a curve of compact type $C=Y_1 \cup \cdots \cup Y_s$ of arithmetic genus $g$ with nodes at the points $\{p_{ij}\}_{ij}$, let $\{l_{Y_1},\dots l_{Y_s} \}$ be a limit linear series $\g^r_d$ on $C$. Let $\{q_{ik} \}_k$ be smooth points on $Y_i$, $i=1,\dots,s$. In \cite{MR846932} a moduli space of such limit series is constructed as a disjoint union of schemes on which the vanishing sequences of the aspects $l_{Y_i}$'s at the nodes are specified. A key property is the additivity of the adjusted Brill-Noether number, that is 
\[
\rho(g,r,d,\{\alpha^{l_{Y_i}}(q_{ik})\}_{ik})\geq \sum_i \rho(Y_i, \{p_{ij}\}_j,\{q_{ik}\}_k).
\]

The smoothing result \cite[Cor. 3.7]{MR846932} assures the smoothability of dimensionally proper limit series. The following facts ease the computations. The adjusted Brill-Noether number for any $\g^r_d$ on one-pointed elliptic curves or on $n$-pointed rational curves is non\-ne\-gative. For a general curve $C$ of arbitrary genus $g$, one has $\rho(C,p)\geq 0$ for $p$ general in $C$ and $\rho(C,y)\geq -1$ for any $y\in C$ (see \cite{MR985853}).

\section{Ramifications on some families of linear series with ${\rho}=0$ or ${\rho}=-1$}
\label{secram}

Here we prove Thm. \ref{ramif0}. The result will be repeatedly used in the next section.

\begin{proof}[Proof of Thm. \ref{ramif0}]
Clearly it is enough to prove the statement for $i=r-1$. We proceed by contradiction. Suppose that for $(C,y)$ a general pointed curve of genus $g$, there exists $x\in C$ such that $h^0(l(- a_{r-1} y-x))\geq 2$, for some $l$ a $\g^r_d$ with $\rho(C,y)=0$. Let us degenerate $C$ to a transversal union $C_1\cup_{y_1} E_1$, where $C_1$ has genus $g-1$ and $E_1$ is an elliptic curve. Since $y$ is a general point, we can assume $y\in E_1$ and $y-y_1$ not to be a $d!$-torsion point in ${\rm Pic}^0(E_1)$. Let $\{ l_{C_1}, l_{E_1} \}$ be a limit $\g^r_d$ on $C_1\cup_{y_1} E_1$ such that $a^{l_{E_1}}(y)=(a_0, a_1, \dots, a_r)$. Denote by $(\alpha_0, \dots, \alpha_r)$ the corresponding ramification sequence. We have that $\rho(C_1,y_1)=\rho(E_1,y,y_1)=0$, hence $w^{l_{C_1}}(y_1)=r +\rho$, where $\rho=\rho(g,r,d)$. Denote by $(b^1_0, b^1_1, \dots, b^1_r)$ the vanishing sequence of $l_{C_1}$ at $y_1$ and by $(\beta^1_0, \beta^1_1, \dots, \beta^1_r)$ the corresponding ramification sequence. 

Suppose $x$ specializes to $E_1$. Then $b^1_r\geq a_r+1$, $b^1_{r-1}\geq a_{r-1}+1$ and we cannot have both equalities, since $y-y_1$ is not in ${\rm Pic}^0(E_1)[d!]$ (see for instance \cite[Prop. 4.1]{MR1785575}). Moreover, as usually $b^1_{k}\geq a_k $ for $0\leq k \leq r-2$, and again among these inequalities there cannot be more than one equality. We deduce
\[
w^{l_{C_1}}(y_1)\geq w^{l_{E_1}}(y)+3+r-2>w^{l_{E_1}}(y)+r=r+\rho 
\]
hence a contradiction. We have supposed that $h^0(l(- a_{r-1} y-x))\geq 2$. Then this pencil degenerates to $l_{E_1}(-a_{r-1}y)$ and to a compatible sub-pencil $l'_{C_1}$ of $l_{C_1}(-x)$. We claim that 
\[
h^0 \left(l_{C_1} \left(-b^1_{r-1} y_1-x \right) \right)\geq 2. 
\]
Suppose this is not the case. Then we have $a^{l_{C_1}(-x)}(y_1) \leq (b^1_0,\dots, b^1_{r-2},b^1_r)$, hence $b^1_r\geq a_r$, $b^1_{r-2}\geq a_{r-1}$ and $b^1_k\geq a_k$, for $0\leq k \leq r-3$. Among these, we cannot have more than one equality, plus $\beta^1_{r-2}\geq \alpha_{r-1}+1$ and $\beta^1_{r-1}\geq \beta^1_{r-2}>\alpha_{r-1}\geq\alpha_{r-2}$, hence 
\[
w^{l_{C_1}}(y_1)\geq w^{l_{E_1}}(y)+1+r-1+\beta^1_{r-1}-\alpha_{r-2}>r+\rho
\]
a contradiction.

From our assumptions, we have deduced that for $(C_1,y_1)$  a general pointed curve of genus $g-1$, there exist $l_{C_1}$ a $\g^r_d$ and $x\in C_1$ such that $\rho(C_1,y_1)=0$ and $h^0(l_{C_1}(-b^1_{r-1} y_1-x))\geq 2$, where $b^1_{r-1}$ is as before.

Then we apply the following recursive argument. At the step $i$, we degenerate the pointed curve $(C_i, y_i)$ of genus $g-i$ to a transversal union $C_{i+1}\cup_{y_{i+1}} E_{i+1}$, where $C_{i+1}$ is a curve of genus $g-i-1$ and $E_{i+1}$ is an elliptic curve, such that $y_i \in E_{i+1}$. Let $\{ l_{C_{i+1}}, l_{E_{i+1}} \}$ be a limit $\g^r_d$ on $C_{i+1}\cup_{y_{i+1}} E_{i+1}$ such that $a^{l_{E_{i+1}}}(y_i)=(b^i_0, b^i_1, \dots, b^i_r)$. From $\rho(C_{i+1},y_{i+1})=\rho(E_{i+1},y_i,y_{i+1})=0$, we compute that $w^{l_{C_{i+1}}}(y_{i+1})=(i+1)r + \rho$. Denote by $(b^{i+1}_0, b^{i+1}_1, \dots, b^{i+1}_r)$ the vanishing sequence of $l_{C_{i+1}}$ at $y_{i+1}$. As before we arrive to a contradiction if $x\in E_{i+1}$, and we deduce
\[
h^0\left(l_{C_{i+1}}\left(- b^{i+1}_{r-1} y_{i+1}-x\right)\right)\geq 2.
\]

At the step $g-2$, our degeneration produces two elliptic curves $C_{g-1}\cup_{y_{g-1}} E_{g-1}$, with $y_{g-2} \in E_{g-1}$. Our assumptions yield the existence of $x \in C_{g-1}$ such that 
\[
h^0 (l_{C_{g-1}}(- b^{g-1}_{r-1} y_{g-1}-x ) )\geq 2.
\]
We compute $w^{l_{C_{i+1}}}(y_{g-1})=(g-1)r + \rho$. By the numerical hypothesis, we see that $(g-1)r + \rho=(d-r-1)(r+1)+1$, hence the vanishing sequence of $l_{C_{g-1}}$ at $y_{g-1}$ has to be $(d-r-1,\dots,d-3,d-2,d)$. Whence the contradiction.
\end{proof}

\noindent The following proves the similar result for some families of linear series with Brill-Noether number $-1$.

\begin{proof}[Proof of Thm \ref{bpf-1}]
The statement says that for every $y\in C$ such that $\rho(C,y)=-1$ for some $l$ a $\g^r_d$, and for every $x\in C$, we have that $h^0(l(- a_{1} y-x))\leq r-1$. This is a closed condition and, using the irreducibility of the divisor $\DD$ of pointed curves admitting a linear series $\g^r_d$ with adjusted Brill-Noether number $-1$, it is enough to prove it for $[C,y]$ general in $\DD$. 

We proceed by contradiction. Suppose for $[C,y]$ general in $\DD$ there exists $x\in C$ such that $h^0(l(- a_{1} y-x))\geq r$ for some $l$ a $\g^r_d$ with $\rho(C,y)=-1$. Let us degenerate $C$ to a transversal union $C_1\cup_{y_1} E_1$ where $C_1$ is a general curve of genus $g-1$ and $E_1$ is an elliptic curve. Since $y$ is a general point, we can assume $y\in E_1$. Let $\{ l_{C_1}, l_{E_1} \}$ be a limit $\g^r_d$ on $C_1\cup_{y_1} E_1$ such that $a^{l_{E_1}}(y)=(a_0, a_1, \dots, a_r)$. Then $\rho(E_1,y,y_1)\leq -1$ and $\rho(C_1,y_1)=0$, hence $w^{l_{C_1}}(y_1)=r+\rho$ (see also \cite[Proof of Thm. 4.6]{MR2530855}). Let $(b^1_0, b^1_1, \dots, b^1_r)$ be the vanishing sequence of $l_{C_1}$ at $y_1$ and $(\beta^1_0, \beta^1_1, \dots, \beta^1_r)$ the corresponding ramification sequence.

The point $x$ has to specialize to $C_1$. Indeed suppose $x\in E_1$. Then $b^1_k \geq a_k +1$ for $k \geq 1$. This implies $w^{l_{C_1}}(y_1)\geq w^{l_{E_1}}(y) +r > \rho +r$, hence a contradiction. Then $x\in C_1$, and $l(- a_{1} y-x)$ degenerates to $l_{E_1}(-a_{1}y)$ and to a compatible system $l'_{C_1} := l_{C_1}(-x)$. We claim that 
\[
h^0 \left( l_{C_1} \left( -b^1_{r-1} y_1-x \right) \right) \geq 2.
\] 
Suppose this is not the case. Then we have $a^{l'_{C_1}}(y_1) \leq (b^1_0,\dots, b^1_{r-2},b^1_r)$ and so $b^1_r\geq a_r$, and $b^1_k\geq a_{k+1}$ for $0\leq k\leq r-2$. Then $\beta^1_k\geq \alpha_{k+1}+1$ for $k\leq r-2$, and summing up we receive
\[
w^{l_{C_1}}(y_1)\geq w^{l_{E_1}}(y)+r-1+\beta^1_{r-1}-\alpha_0. 
\]
Clearly $\beta^1_{r-1}\geq \beta^1_{r-2}>\alpha_{r-1}\geq \alpha_0$. Hence $w^{l_{C_1}}(y_1) > \rho +r$, a contradiction.

All in all from our assumptions we have deduced that for a general pointed curve $(C_1,y_1)$ of genus $g-1$, there exist $l_{C_1}$ a $\g^r_d$ and $x\in C_1$ such that $\rho(C_1,y_1)=0$ and $h^0(l_{C_1}(-b^1_{r-1} y_1-x))\geq 2$, where $b^1_{r-1}$ is as before. This contradicts Thm. \ref{ramif0}, hence we receive the statement.
\end{proof}

\section{The divisor $\D^2_d$}
\label{class}

\noindent Remember that $\mbox{Pic}_\QQ(\Mbar_{g,1})$ is generated by the Hodge class $\lambda$, the cotangent class $\psi$ corresponding to the marked point, and the boundary classes $\delta_0,\dots \delta_{g-1}$ defined as follows. The class $\delta_0$ is the class of the closure of the locus of pointed irreducible nodal curves, and the class $\delta_i$ is the class of the closure of the locus of pointed curves $[C_i\cup C_{g-i},p]$ where $C_i$ and $C_{g-i}$ are smooth curves respectively of genus $i$ and $g-i$ meeting transversally in one point, and $p$ is a smooth point in $C_i$, for $i=1,\dots,g-1$. In this section we prove the following theorem.

\begin{thm}
Let $g=3s$ and $d=2s+2$ for $s\geq1$. The class of the divisor $\overline{\D}^2_d$ in $\mbox{\rm Pic}_\QQ(\Mbar_{g,1})$ is
\[
\left[\overline{ \D}^2_d \right] = a \lambda + c \psi - \sum_{i=0}^{g-1} b_i \delta_i 
\]
where
\begin{eqnarray*}
a &=& \frac{48 s^4 + 80 s^3 -16 s^2 - 64s +24}{(3s-1)(3s-2)(s+3)} N_{g,2,d} \\
c & = & \frac{2s(s-1)}{3s-1} N_{g,2,d}\\
b_0 &=& \frac{24 s^4 + 23 s^3 - 18 s^2 - 11s +6}{3(3s-1)(3s-2)(s+3)} N_{g,2,d}\\
b_1 & = & \frac{14 s^3 + 6 s^2 - 8 s}{(3s-2)(s+3)} N_{g,2,d} \\
b_{g-1} & = & \frac{48 s^4 + 12 s^3 - 56 s^2 +20s}{(3s-1)(3s-2)(s+3)} N_{g,2,d}.
\end{eqnarray*}
Moreover for $g=6$ and for $i=2,3,4,$ we have that
\begin{eqnarray*}
b_i = -7 i^2+43 i-6.
\end{eqnarray*}
\end{thm}

\subsection{The coefficient $c$} 

The coefficient $c$ can be quickly found. Let $C$ be a general curve of genus $g$ and consider the curve $\overline{C} = \{[C,y] : y \in C\}$ in $\Mbar_{g,1}$ obtained varying the point $y$ on $C$. Then the only generator class having non-zero intersection with $\overline{C}$ is $\psi$, and $\overline{C} \cdot \psi = 2g-2$. On the other hand, $\overline{C} \cdot \overline{\D}^2_d$ is equal to the number of triples $(x,y,l)\in C \times C \times G^2_d(C)$ such that $x$ and $y$ are different points and $h^0(l(-x-y))\geq 2$. The number of such linear series on a general $C$ is computed by the Castelnuovo number (remember that $\rho =0$), and for each of them the number of couples $(x,y)$ imposing only one condition is twice the number of double points, computed by the Pl\"ucker formula. Hence we get the equation
\[
\overline{\D}^2_d \cdot \overline{C} =  2 \left( \frac{(d-1)(d-2)}{2} -g \right) N_{g,2,d} = c \, (2g-2)
\]
\noindent and so 
\[
c= \frac{2s(s-1)}{3s-1} N_{g,2,d}.
\]

\subsection{The coefficients $a$ and $b_0$} 

In order to compute $a$ and $b_0$, we use a Porteous-style argument. Let $\GG^2_d$ be the family parametrizing triples $(C,p,l)$, where $[C,p]\in \M_{g,1}^{\rm irr}$ and $l$ is a $\g^2_d$ on $C$; denote by $\eta: \GG^2_d \rightarrow \M_{g,1}^{\rm irr}$ the natural map. There exists $\pi:\YY^2_d \rightarrow \GG^2_d$ a universal pointed quasi-stable curve, with $\sigma: \GG^2_d \rightarrow \YY^2_d$ the marked section. Let $\Li\rightarrow \YY^2_d$ be the universal line bundle of relative degree $d$ together with the trivialization $\sigma^*(\Li) \cong \OO_{ \GG^2_d}$, and $\VV \subset \pi_*(\Li)$ be the sub-bundle which over each point $(C,p,l=(L,V))$ in $\GG^2_d$ restricts to $V$.
(See \cite[\S 2]{khosla} for more details.)
 
Furthermore let us denote by $\Zz^2_d$ the family parametrizing $\left( (C,p), x_1, x_2, l  \right)$, where $[C,p]\in \M^{\rm irr}_{g,1}$, $x_1,x_2 \in C$ and $l$ is a $\g^2_d$ on $C$, and let $\mu , \nu : \Zz^2_d \rightarrow \YY^2_d$ be defined as the maps that send $\left( (C,p), x_1, x_2, l  \right)$ respectively to $\left( (C,p), x_1,  l  \right)$ and $\left( (C,p), x_2, l  \right)$.

Now given a linear series $l=(L,V)$, the natural map
\[
\varphi : V \rightarrow H^0(L |_{p+x})
\]
globalizes to
\[
\widetilde{\varphi}: \VV \rightarrow \mu_*\left( \nu^* \Li \otimes \OO/ \II_{ \Gamma_\sigma + \Delta}  \right) =: \MM
\]
as a map of vector bundle over $\YY^2_d$, where $\Delta$ and $\Gamma_\sigma$ are the loci in $\Zz^2_d$ determined respectively by $x_1=x_2$ and $x_2 = p$. Then $\overline{\D}^2_d \cap \M^{\rm irr}_{g,1}$ is the push-forward of the locus in $\YY^2_d$ where $\widetilde{\varphi}$ has rank $\leq 1$. Using Porteous formula, we have
\begin{eqnarray}
\label{class D}
 [\overline{\D}^2_d] | _{\M^{\rm irr}_{g,1}} & = & \eta_* \pi_* \left[  \frac{\VV^\vee}{\MM^\vee}  \right]_2  \\
					    & = & \eta_* \pi_*  \left( \pi^*c_2(\VV^\vee) + \pi^*c_1(\VV^\vee) \cdot c_1(\MM) + c_1^2(\MM) - c_2(\MM) \right)  .\nonumber 
\end{eqnarray}
Let us find the Chern classes of $\MM$. Tensoring the exact sequence
\[
0 \rightarrow \II_\Delta/ \II_{\Delta + \Gamma_\sigma} \rightarrow \OO / \II_{\Delta + \Gamma_\sigma} \rightarrow \OO_\Delta \rightarrow 0
\]
\noindent by $\nu^* {\Li}$ and applying $\mu_*$, we deduce that
\begin{eqnarray*}
ch(\MM) &=& ch(\mu_*(\OO_{\Gamma_\sigma}(-\Delta) \otimes \nu^*\Li)) + ch(\mu_*(\OO_\Delta \otimes \nu^*\Li))  \\
	      &=& ch(\mu_*(\OO_{\Gamma_\sigma}(-\Delta))) + ch(\mu_*(\OO_\Delta \otimes \nu^*\Li))  \\
              &=& e^{-\sigma} +  ch(\Li) 
\end{eqnarray*}
hence
\begin{eqnarray*}
c_1(\MM)  &  =  &  c_1(\Li) - \sigma  \\
c_2(\MM)  &  =  &  -\sigma c_1(\Li).
\end{eqnarray*}
The following classes 
\begin{eqnarray*}
\alpha & = & \pi_* \left( c_1(\Li)^2 \cap [\YY^2_d]  \right) \\
\gamma & = & c_1(\VV) \cap [\GG^2_d]
\end{eqnarray*}
have been studied in \cite[Thm. 2.11]{khosla}. In particular
\begin{eqnarray*}
\frac{6(g-1)(g-2)}{d N_{g,2,d}} \, \eta_*(\alpha) | _{\M^{\rm irr}_{g,1}} & = &  6(gd-2g^2+8d-8g+4) \lambda \\
								    &   & {}+ (2g^2-gd+3g-4d-2)\delta_0\\
								    &   & {}-6d(g-2) \psi ,\\
\frac{2(g-1)(g-2)}{N_{g,2,d}} \, \eta_*(\gamma) | _{\M^{\rm irr}_{g,1}} & = &   \left( -(g+3)\xi +40  \right) \lambda \\
								  &   &   {}+\frac{1}{6} \left( (g+1)\xi -24  \right) \delta_0 \\
								  &   &   {}-3d(g-2)\psi , 
\end{eqnarray*}
where
\[
\xi = 3(g-1) + \frac{(g+3)(3g-2d-1)}{g-d+5}.
\]
Plugging into (\ref{class D}) and using the projection formula, we find
\begin{eqnarray*}
[\overline{\D}^2_d] | _{\M^{\rm irr}_{g,1}}  &=& \eta_* \left( -\gamma \cdot \pi_*c_1(\Li)  + \gamma \cdot \pi_*\sigma + \alpha +\pi_*\sigma^2 - \pi_*(\sigma c_1(\Li))   \right) \\
				  &=& (1-d) \eta_* (\gamma) + \eta_* (\alpha) -N_{g,2,d}\cdot \psi .
\end{eqnarray*}
Hence
\begin{eqnarray*}
a &=& \frac{48 s^4 + 80 s^3 -16 s^2 - 64s +24}{(3s-1)(3s-2)(s+3)} N_{g,2,d} \\
b_0 &=& \frac{24 s^4 + 23 s^3 - 18 s^2 - 11s +6}{3(3s-1)(3s-2)(s+3)} N_{g,2,d}
\end{eqnarray*}
and we recover the previously computed coefficient $c$.

\subsection{The coefficient $b_1$} 
\label{b_1}

Let $C$ be a general curve of genus $g-1$ and $(E,p,q)$ a two-pointed elliptic curve, with $p-q$ not a torsion point in Pic$^0(E)$. Let $\overline{C}_1:=\{ (C \cup_{y \sim q} E, p) \}_{y\in C}$ be the family of curves obtained identifying the point $q\in E$ with a moving point $y\in C$. Computing the intersection of the divisor $\overline{\D}^2_d$ with $\overline{C}_1$ is equivalent to answering the following question: how many triples $(x,y,l)$ are there, with $y\in C$, $x\in C \cup_{y \sim q} E \setminus \{ p\}$ and $l=\{l_C,l_E\}$ a limit $\g^2_d$ on $C \cup_{y \sim q} E$, such that $(p,x,l)$ arises as limit of $(p_t, x_t, l_t)$ on a family of curves $\{C_t\}_t$ with smooth general element, where $p_t$ and $x_t$ impose only one condition on $l_t$ a $\g^2_d$?

\begin{center}
\includegraphics[width=4.5cm]{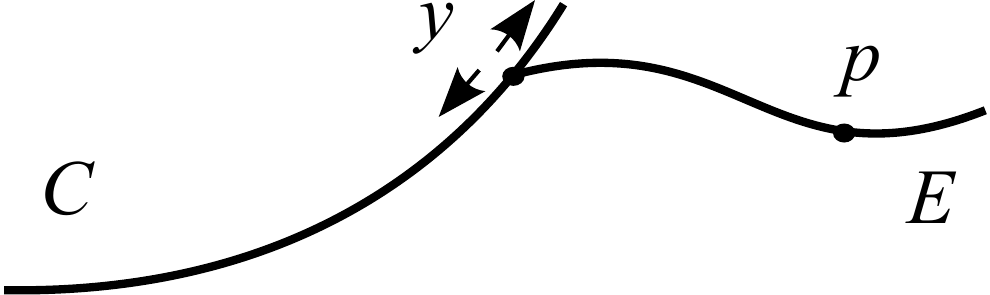}
\end{center}

Let $a^{l_E}(q)=(a_0, a_1, a_2)$ be the vanishing sequence of $l_E \in G^2_d(E)$ at $q$. Since $C$ is general, there are no $\g^2_{d-1}$ on $C$, hence $l_C$ is base-point free and $a_2=d$. Moreover we know $a_1\leq d-2$.
Let us suppose $x\in E\setminus \{ q\}$. We distinguish two cases. If $\rho(E,q)= \rho(C,y) = 0$, then $w^{l_E}(q)=  \rho(1,2,d)=3d-8$.  Thus $a^{l_E}(q)=(d-3,d-2,d)$. Removing the base point we have that $l_E (- (d-3)q)$ is a $\g^2_3$ and $l_E (-(d-3)q-p-x)$ produces a $\g^1_1$ on $E$, hence a contradiction. The other case is $\rho(E, q)=1$ and $\rho(C,y)\leq-1$. These force $a^{l_E}(q)=(d-4,d-2,d)$ and $a^{l_C}(y)\geq(0,2,4)$. On $E$ we have that $l_E (-(d-4)q-p-x)$ is a $\g^1_2$. 

The question splits in two: firstly, how many linear series $l_E \in G^2_4(E)$ and points $x\in E\setminus \{q\}$ are there such that $a^{l_E}(q) = (0,2,4)$ and $l_E (-p-x) \in G^1_2(E)$? The first condition restricts our attention to the linear series $l_E=(\OO(4q),V)$ where $V$ is a tridimensional vector space and $H^0(\OO(4q-2q))\subset V$, while the second condition tells us $H^0(\OO(4q-p-x))\subset V$.
If $x=p$, then we get $p-q$ is a torsion point in Pic$^0(E)$, a contradiction. On the other hand, if $x\in E\setminus \{p,q\}$, then $H^0(\OO(4q-2q)) \cap H^0(\OO(4q-p-x)) \not= \emptyset$ entails $p+x \equiv 2q$. Hence the point $x$ and the space $V=H^0(\OO(4q-2q)) + H^0(\OO(4q-p-x))$ are uniquely determined.

Secondly, how many couples $(y,l_C)\in C\times G^2_d(C)$ are there, such that the va\-ni\-shing sequence of $l_C$ at $y$ is greater than or equal to $(0,2,4)$? This is a particular case of a problem discussed in \cite[Proof of Thm. 4.6]{MR2530855}. The answer is
\begin{multline*}
(g-1) \left( 15N_{g-1,2,d,(0,2,2)} + 3N_{g-1,2,d,(1,1,2)} + 3N_{g-1,2,d,(0,1,3)}  \right) \\= \frac{24(2s^2+3s-4)}{s+3} N_{g,2,d}.
\end{multline*}

Now let us suppose $x\in C\setminus \{ y\}$. The condition on $x$ and $p$ can be reformulated in the following manner. We consider the curve $C\cup_y E$ as the special fiber $X_0$ of a family of curves $\pi:X \rightarrow B$ with sections $x(t)$ and $p(t)$ such that $x(0)=x$, $p(0)=p$, and with smooth general fiber having $l=(\Li, V)$ a $\g^2_d$ such that $l(-x-p)$ is a $\g^1_{d-2}$. Let $V'\subset V$ be the two dimensional linear subspace formed by those sections $\sigma\in V$ such that div$(\sigma)\geq x+p$. Then $V'$ specializes on $X_0$ to $V'_C\subset V_C$ and $V'_E \subset V_E$ two-dimensional subspaces, where $\{ l_{C}=(\Li_C, V_C), l_{E}=(\Li_E, V_E)\}$ is a limit $\g^2_d$, such that 
\[\left\{
\begin{array}{l}
{\rm ord}_y(\sigma_C) + {\rm ord}_y(\sigma_E) \geq d \\
{\rm div}(\sigma_C) \geq x \\
{\rm div}(\sigma_E) \geq p
\end{array}
\right.\]

\noindent for every $\sigma_C\in V'_C$ and $\sigma_E\in V'_E$. Let $l'_C:=(\Li_C, V'_C)$ and $l'_E:=(\Li_E, V'_E)$. Note that since $\sigma_E \geq p$, we get ${ \rm ord}_y(\sigma_E)<d$, $\forall\, \sigma_E\in V'_E$.  Then ${\rm ord}_y (\sigma_C)>0$, hence ${\rm ord}_y (\sigma_C)\geq 2$, since $y$ is a cuspidal point on $C$. Removing the base point,  $l'_C$ is a $\g^1_{d-2}$ such that $l'_C (-x)$ is a $\g^1_{d-3}$. Let us suppose $\rho(E,y)=1$ and $\rho(C,y)=-1$. Then $a^{l_{E}}(y)=(d-4,d-2,d)$, $a^{l'_{E}}(y)=(d-4,d-2)$, $a^{l_{C}}(y)=(0,2,4)$ and $a^{l'_{C}}(y)=(2,4)$. Now $l_{C}$ is characterized by the conditions $H^0(l_{C}(-2y-x))\geq 2$ and $H^0(l_{C}(-4y-x))\geq 1$. By Thm. $\ref{bpf-1}$ this possibility does not occur.

Suppose now $\rho(E,y)=\rho(C,y)=0$. Then $a^{l_{E}}(y)=(d-3,d-2,d)$, i.e. $l_E(-(d-3)y)=|3y|$ is uniquely determined. On the $C$ aspect we have that $a^{l_{C}}(y)=(0,2,3)$ and $h^0(l_C(-2y-x))\geq 2$. Hence we are interested on $Y$ the locus of triples $(x,y,l_C)$ such that the map
\[
\varphi : H^0(l_C) \rightarrow H^0(l_C|_{2y+x})
\]
\noindent has rank $\leq 1$. By Thm. \ref{ramif0} there is only a finite number of such triples, and clearly the case $a^{l_C}(y)>(0,2,3)$ cannot occur. Moreover, note that $x$ and $y$ will be necessarily distinct. 

Let $\mu = \pi_{1,2,4}: C\times C \times C \times W^2_d(C) \rightarrow C\times C \times W^2_d(C)$ and $\nu = \pi_{3,4}: C\times C \times C \times W^2_d(C) \rightarrow C \times W^2_d(C)$ be the natural projections respectively on the first, second and forth components, and on the third and forth components. Let $\pi: C\times C \times W^2_d(C) \rightarrow W^2_d(C)$ be the natural projection on the third component.  Now $\varphi$ globalizes to
\[
\widetilde{\varphi}: \pi^* \EE \rightarrow \mu_* \left( \nu^* \Li \otimes \OO/\II_\DD \right) =: \MM
\]
\noindent as a map of rank $3$ bundles over $C\times C \times W^2_d(C)$, where $\DD$ is the pullback to $C\times C\times C \times W^2_d(C)$ of the divisor on $C\times C \times C$ that on $(x,y,C)\cong C$ restricts to $x+2y$, $\Li$ is a Poincar${\rm \acute{e}}$ bundle on $C\times W^2_d$ and $\EE$ is the push-forward of $\Li$ to $W^2_d(C)$. Then $Y$ is the degeneracy locus where $\widetilde{\varphi}$ has rank $\leq 1$.
Let $\cc_i:= c_i(\EE)$ be the Chern classes of $\EE$. By Porteous formula, we have
\[
[Y]= \left[
\begin{array}{cc}
e_2 & e_3 \\
e_1 & e_2
\end{array}
\right]
\]
\noindent where the $e_i$'s are the Chern classes of $ \pi^* \EE^\vee - \MM^\vee$,  i.e.
\begin{eqnarray*}
e_{1} & = & \cc_{1} + c_{1}(\MM) \nonumber \\
e_{2} & = & \cc_{2} + \cc_{1}c_{1}(\MM) + c_{1}^2(\MM) - c_{2}(\MM) \\ 
e_{3} & = & \cc_{3} + \cc_{2}c_{1}(\MM) + \cc_{1}\left( c_{1}^2(\MM)-c_{2}(\MM) \right)  \nonumber \\ 
        	 &    &{}+\left( c_{1}^3(\MM) + c_{3}(\MM) -2c_{1}(\MM)c_{2}(\MM) \right). \nonumber
\end{eqnarray*}

Let us find the Chern classes of $\MM$. First we develop some notations (see also \cite[\S VIII.2]{MR770932}). Let $\pi_i:C\times C\times C\times W^2_d(C) \rightarrow C$ for $i=1,2,3$ and $\pi_4:C\times C\times C\times W^2_d(C) \rightarrow W^2_d(C)$ be the natural projections. Denote by $\theta$ the pull-back to $C\times C\times C\times W^2_d(C)$ of the class $\theta \in H^2(W^2_d(C))$ via $\pi_4$, and denote by $\eta_i$ the cohomology class $\pi_i^*([{\rm point}]) \in H^2(C\times C\times C\times W^2_d(C))$, for $i=1,2,3$. Note that $\eta_i^2=0$. Furthermore, given a symplectic basis $\delta_1,\dots,\delta_{2(g-1)}$ for $H^1(C,\ZZ)\cong H^1(W^2_d(C),\ZZ)$, denote by $\delta^i_\alpha$ the pull-back to $C\times C\times C\times W^2_d(C)$ of $\delta_\alpha$ via $\pi_i$, for $i=1,2,3,4$. Let us define
\[
\gamma_{ij}:= - \sum_{\alpha=1}^{g-1} \left( \delta^j_\alpha \delta^i_{g-1+\alpha} -  \delta^j_{g-1+\alpha} \delta^i_\alpha   \right).
\]
Note that 
\[
\begin{array}{cclccccclc}
\gamma_{ij}^2 &=& -2(g-1) \eta_i \eta_j  & \mbox{ and}  & \eta_i \gamma_{ij} &=& \gamma_{ij}^3=0   &\mbox{ for}  & 1\leq i<j \leq 3,\\
\gamma_{k4}^2 &=& -2 \eta_k \theta & \mbox{ and} & \eta_k \gamma_{k4} &=& \gamma_{k4}^3=0 & \mbox{ for} & k=1,2,3.
\end{array}
\]
Moreover
\[
\gamma_{ij}\gamma_{jk}  \quad = \quad  \eta_j \gamma_{ik}, 
\]
for $1\leq i < j < k \leq 4$. With these notations, we have
\[
ch(\nu ^* \Li \otimes \OO / \II_\DD) = (1+d\eta_3 + \gamma_{34}-\eta_3 \theta)\left(1-e^{-(\eta_1+\gamma_{13}+\eta_3 + 2\eta_2 + 2\gamma_{23}+2\eta_3)}\right),
\]
\noindent hence by Grothendieck-Riemann-Roch
\begin{eqnarray*}
ch(\MM) &=& \mu_* \left(  (1+(2-g)\eta_3) ch(\nu ^* \Li \otimes \OO / \II_\DD) \right) \\
              &=& 3+(d-2)\eta_1 + (2g+2d-6)\eta_2 -2 \gamma_{12} + \gamma_{14} + 2 \gamma_{24} \\
              && {}-\eta_1 \theta - 2\eta_2 \theta + (8-2d-4g)\eta_1 \eta_2 - 2 \eta_1 \gamma_{24} -2\eta_2 \gamma_{14} + 2\eta_1 \eta_2 \theta.
\end{eqnarray*}
Using Newton's identities, we recover the Chern classes of $\MM$:
\begin{eqnarray*}
c_1(\MM) &=& (d-2)\eta_1 + (2g+2d-6)\eta_2 - 2 \gamma_{12} + \gamma_{14} + 2\gamma_{24}, \\
c_2(\MM) &=& (2d^2-8d +2gd+8 -4g )\eta_1\eta_2 + (2g+2d-8)\eta_2\gamma_{14} \\
                 && {}+ (2d-4)\eta_1\gamma_{24} +2\gamma_{14}\gamma_{24} -2\eta_2\theta, \\
c_3(\MM) &=& (4- 2d)\eta_1\eta_2\theta  - 2\eta_2\gamma_{14}\theta.
\end{eqnarray*}
We finally find
\begin{eqnarray*}
[Y] &=& \eta_1 \eta_2  (\cc_1^2 (2d^2 -8d +2dg+4-4(g-1) )  \\ 
      && {}+\cc_1 \theta (-12d -4g+40) +\cc_2(-4d+16-8g) +12 \theta^2 )\\
      &=& \frac{(28s+48)(s-2)(s-1)}{(s+3)}N_{g,2,d} \cdot  \eta_1 \eta_2  \theta^{g-1}
\end{eqnarray*}
where we have used the following identities proved in \cite[Lemma 2.6]{MR2530855}
\begin{eqnarray*}
\cc_1^2 & = & \left(1+\frac{2s+2}{s+3}\right)\cc_2 \\
\cc_1 \theta & = & (s+1)\cc_2\\
\theta^2 &=& \frac{(s+1)(s+2)}{3} \cc_2\\
\cc_2 &=& N_{g,2,d} \cdot \theta^{g-1}.
\end{eqnarray*}

We are going to show that we have already considered all non zero contributions. 
Indeed let us suppose $x=y$. Blowing up the point $x$, we obtain $C\cup_y\PP^1\cup_q E$ with $x\in \PP^1\setminus \{y,q\}$ and $p\in E\setminus \{q\}$. We reformulate the condition on $x$ and $p$ viewing our curve as the special fiber of a family of curves $\pi:X\rightarrow B$ as before. Let $\{l_C, l_{\PP^1}, l_E\}$ be a limit $\g^2_d$. Now $V'$ specializes to $V'_C$, $V'_{\PP^1}$ and $V'_E$. There are three possibilities: either $\rho(C,y)=\rho(\PP^1,x,y,q)=\rho(E,p,q)=0$, or $\rho(C,y)=-1$, $\rho(\PP^1,x,y,q)=0$, $\rho(E,p,q)=1$, or $\rho(C,y)=-1$, $\rho(\PP^1,x,y,q)=1$, $\rho(E,p,q)=0$. In all these cases $a^{l_C}(y)=(0,2,a_2^{l_C}(y))$ (remember that $l_C$ is base-point free) and $a^{l_E}(q)=(a_0^{l_E}(q),d-2,d)$. Hence $a^{l_{\PP^1}}(y)=(a_0^{l_{\PP^1}}(y),d-2,d)$ and $a^{l_{\PP^1}}(q)=(0,2,a_2^{l_{\PP^1}}(q))$. Let us restrict now to the sections in $V'_C$, $V'_{\PP^1}$ and $V'_E$. For all sections $\sigma_{\PP^1}\in V'_{\PP^1}$ since ${\rm div}(\sigma_{\PP^1})\geq x$, we have that ${\rm ord}_y (\sigma_{\PP^1})<d$ and hence ${\rm ord}_y(\sigma_{\PP^1})\leq d-2$. On the other side, since for all $\sigma_E\in V'_E$, ${\rm div}(\sigma_E)\geq p$, we have that ${\rm ord}_q (\sigma_E)<d$ and hence ${\rm ord}_q(\sigma_{\PP^1})\geq 2$. Let us take one section $\tau\in V'_{\PP^1}$ such that ${\rm ord}_y(\tau)=d-2$. Since ${\rm div}(\tau)\geq (d-2)y+x$, we get ${\rm ord}_q(\tau)\leq 1$, hence a contradiction.

Thus we have that
\[
\overline{\D}^2_d \cdot \overline{C}_1 = \frac{24(2s^2+3s-4)}{s+3}N_{g,2,d} + \frac{(28s+48)(s-2)(s-1)}{(s+3)}N_{g,2,d}.
\]
while considering the intersection of the test curve $\overline{C}_1$ with the ge\-ne\-ra\-ting classes we have
\[
\overline{\D}^2_d \cdot \overline{C}_1 = b_1 (2g-4),
\]
whence
\[
b_1 = \frac{14 s^3 + 6 s^2 - 8 s}{(3s-2)(s+3)} N_{g,2,d} .
\]

\begin{rem}
The previous class $[Y]$ being nonzero, it implies together with Thm. \ref{ramif0} that the scheme $\GG^2_d((0,2,3))$ over $\M_{g-1,1}$ is not irreducible.
\end{rem}

\subsection{The coefficient $b_{g-1}$} 

We analyze now the following test curve $\overline{E}$. Let $(C,p)$ be a general pointed curve of genus $g-1$ and $(E,q)$ be a pointed elliptic curve. Let us identify the points $p$ and $q$ and let $y$ be a movable point in $E$. We have
\[
0 = \overline{\D}^2_d \cdot \overline{E} = c+ b_1 -b_{g-1},
\]
whence
\[
b_{g-1} = \frac{48 s^4 + 12 s^3 - 56 s^2 +20s}{(3s-1)(3s-2)(s+3)} N_{g,2,d}.
\]

\subsection{A test} 

Furthermore, as a test we consider the family of curves $R$. Let $(C,p,q)$ be a general two-pointed curve of genus $g-1$ and let us identify the point $q$ with the base point of a general pencil of plane cubic curves. 
We have
\[
0 = \overline{\D}^2_d \cdot R = a - 12 b_0 + b_{g-1}.
\]

\subsection{The remaining coefficients in case $g=6$}

Denote by $P_g$ the moduli space of stable $g$-pointed rational curves. Let $(E,p,q)$ be a general $2$-pointed elliptic curve and let $j: P_g \rightarrow \Mbar_{g,1}$ be the map obtained identifying the first marked point on a rational curve with the point $q\in E$ and attaching a fixed elliptic tail at the other marked points. We claim that $j^*(\overline{\D}^2_6)=0$.

Indeed consider a flag curve of genus $6$ in the image of $j$. Clearly the only possibility for the adjusted Brill-Noether numbers is to be zero on each aspect. In particular the collection of the aspects on all components but $E$ smooths to a $\g^2_6$ on a general $1$-pointed curve of genus $5$.  As discussed in section \ref{b_1}, the point $x$ can not be in $E$. Suppose $x$ is in the rest of the curve. Then smoothing we get $l$ a $\g^2_6$ on a general pointed curve of genus $5$ such that $l(-2q-x))$ is a $\g^1_3$, a contradiction.

Now let us study the pull-back of the generating classes. As in \cite[\S 3]{MR910206} we have that $j^*(\lambda)=j^*(\delta_0)=0$. Furthermore $j^*(\psi)=0$.

For $i=1,\dots,g-3$ denote by $\varepsilon^{(1)}_i$ the class of the divisor which is the closure in $P_g$ of the locus of $2$-component curves having exactly the first marked point and other $i$ marked points on one of the two components. Then clearly $j^*(\delta_i)= \varepsilon^{(1)}_{i-1}$ for $i=2,\dots,g-2$. Moreover adapting the argument in  \cite[pg. 49]{MR985853}, we have that
\[
j^*(\delta_{g-1}) = -\sum_{i=1}^{g-3} \frac{i(g-i-1)}{g-2} \varepsilon^{(1)}_{i}
\]
while 
\[
j^*(\delta_{1}) = -\sum_{i=1}^{g-3} \frac{(g-i-1)(g-i-2)}{(g-1)(g-2)}  \varepsilon^{(1)}_{i}.
\]
Finally since $j^*(\overline{\D}^2_6)=0$, checking the coefficient of $\varepsilon^{(1)}_{i}$ we obtain
\[
b_{i+1} =   \frac{(g-i-1)(g-i-2)}{(g-1)(g-2)} b_1 + \frac{i(g-i-1)}{g-2} b_{g-1} 
\]
for $i=1,2,3$.

\nocite{MR2530855} \nocite{MR791679} \nocite{MR985853} \nocite{MR1785575} \nocite{MR846932} \nocite{MR993171} \nocite{MR1953519} \nocite{MR1016424} \nocite{MR770932}

\bibliographystyle{alpha}
\bibliography{Biblio.bib}

\begin{thebibliography}{ACGH85}

\bibitem[ACGH85]{MR770932}
E.~Arbarello, M.~Cornalba, P.~A. Griffiths, and J.~Harris.
\newblock {\em Geometry of algebraic curves. {V}ol. {I}}, volume 267 of {\em
  Grundlehren der Mathematischen Wissenschaften [Fundamental Principles of
  Mathematical Sciences]}.
\newblock Springer-Verlag, New York, 1985.

\bibitem[Cuk89]{MR1016424}
Fernando Cukierman.
\newblock Families of {W}eierstrass points.
\newblock {\em Duke Math. J.}, 58(2):317--346, 1989.

\bibitem[Dia85]{MR791679}
Steven Diaz.
\newblock Exceptional {W}eierstrass points and the divisor on moduli space that
  they define.
\newblock {\em Mem. Amer. Math. Soc.}, 56(327):iv+69, 1985.

\bibitem[EH86]{MR846932}
David Eisenbud and Joe Harris.
\newblock Limit linear series: basic theory.
\newblock {\em Invent. Math.}, 85(2):337--371, 1986.

\bibitem[EH87]{MR910206}
David Eisenbud and Joe Harris.
\newblock The {K}odaira dimension of the moduli space of curves of genus {$\geq
  23$}.
\newblock {\em Invent. Math.}, 90(2):359--387, 1987.

\bibitem[EH89]{MR985853}
David Eisenbud and Joe Harris.
\newblock Irreducibility of some families of linear series with
  {B}rill-{N}oether number {$-1$}.
\newblock {\em Ann. Sci. \'Ecole Norm. Sup. (4)}, 22(1):33--53, 1989.

\bibitem[Far00]{MR1785575}
Gavril Farkas.
\newblock The geometry of the moduli space of curves of genus 23.
\newblock {\em Math. Ann.}, 318(1):43--65, 2000.

\bibitem[Far09]{MR2530855}
Gavril Farkas.
\newblock Koszul divisors on moduli spaces of curves.
\newblock {\em Amer. J. Math.}, 131(3):819--867, 2009.

\bibitem[Jen10]{Jensen}
David Jensen.
\newblock Rational fibrations of ${\Mbar}_{5,1}$ and ${\Mbar}_{6,1}$.
\newblock {\em Preprint, arXiv:1012.5115}, 2010.

\bibitem[Kho07]{khosla}
Deepak Khosla.
\newblock Tautological classes on moduli spaces of curves with linear series
  and a push-forward formula when $\rho=0$.
\newblock {\em Preprint, arXiv:0704.1340}, 2007.

\bibitem[Log03]{MR1953519}
Adam Logan.
\newblock The {K}odaira dimension of moduli spaces of curves with marked
  points.
\newblock {\em Amer. J. Math.}, 125(1):105--138, 2003.

\bibitem[SB89]{MR993171}
N.~I. Shepherd-Barron.
\newblock Invariant theory for {$S_5$} and the rationality of {$M_6$}.
\newblock {\em Compositio Math.}, 70(1):13--25, 1989.

\end{thebibliography}

\end{document}